\documentclass{amsart}

\usepackage{amsmath, amssymb, latexsym}
\usepackage{epsfig}
\usepackage{dsfont}
\usepackage{bm}
\usepackage{url}
\usepackage{cite}

\newcommand{\HCurlZ}{H_0(\mathbf{curl};\Omega)}
\newcommand{\HDiv}{H(\mathrm{div})}
\newcommand{\HCurl}{H(\mathbf{curl})}

\newcommand{\CURL}{\mathbf{curl}\,}
\newcommand{\bu}{\bm{u}}
\newcommand{\bv}{\bm{v}}
\newcommand{\bw}{\bm{w}}

\newcommand{\bz}{\bm{z}}

\newcommand{\cT}{\mathcal{T}}
\newcommand{\cF}{\mathcal{F}}
\newcommand{\cE}{\mathcal{E}}
\newcommand{\cV}{\mathcal{V}}

\makeatletter
\@namedef{subjclassname@2020}{%
  \textup{2020} Mathematics Subject Classification}
\makeatother

\begin{document}

\title[Smoothers Based on Nonoverlapping DDM for $\HCurl$]
{Smoothers Based on Nonoverlapping Domain Decomposition Methods for $\HCurl$ Problems: A Numerical Study}

\author{Duk-Soon Oh}
\address{Department of Mathematics\\
  Chungnam National University\\South Korea}
\email{duksoon@cnu.ac.kr}
\thanks{This work was supported by the National Research Foundation of Korea (NRF) grant funded by the Korea government (MSIT) (No. 2020R1F1A1A01072168)}

\subjclass[2020]{65N30, 65N55, 65F08}
\keywords{multigrid method, N\'{e}d\'{e}lec finite element, $\HCurl$, nonoverlapping domain decomposition}

\begin{abstract}
	This paper presents a numerical study on multigrid algorithms of $V$--cycle type for problems posed in the Hilbert space $\HCurl$ in three dimensions. The multigrid methods are designed for discrete problems originated from the discretization using the hexahedral N\'{e}d\'{e}lec edge element of the lowest-order. Our suggested methods are associated with smoothers constructed by substructuring based on domain decomposition methods of nonoverlapping type. Numerical experiments to demonstrate the robustness and the effectiveness of the suggested algorithms are also provided.
	\end{abstract}

	\maketitle

	\section{Introduction}\label{sec:introduction}
	Let $\Omega$ be a domain that is bounded in $\mathds{R}^3$. We will work with the $\HCurlZ$ Hilbert space which consists of vector fields in the space $(L^2(\Omega))^3$ with curl also in $(L^2(\Omega))^3$ and vanishing tangential components on the boundary $\partial \Omega$ (cf. \cite{GR:1986:NS}). We will consider the following variational problem: Find $\bu \in \HCurlZ$ such that
	\begin{equation}\label{eq:ModelProblem}
	a(\bu, \bv) = (\bm{f}, \bv), \quad \forall \bv \in \HCurlZ,
	\end{equation}
	where
	\begin{equation}
	a(\bv, \bw) = \alpha \cdot (\CURL \bv, \CURL \bw) + \beta \cdot (\bv, \bw).
	\end{equation}
	Here, $(\cdot,\cdot)$ is the standard inner product on $[L_2(\Omega)]^3$ and we assume that $\alpha$ is nonnegative and $\beta$ is strictly positive. We also assume that $\bm{f}$ is a square integrable vector field on $\Omega$, i.e., $\bm{f} \in (L^2(\Omega))^3$. In this manuscript, we will provide a multigrid framework for solving our model problem (\ref{eq:ModelProblem}).
	
	The model problem (\ref{eq:ModelProblem}) is originated from the applications in Maxwell's equation; see \cite{Monk:1991:Maxwell}. Relevant fast solvers, such as multigrid methods and domain decomposition methods, for problems connected with $\HCurl$ have been discussed in \cite{Hiptmair:1999:MGHCurl, HX:2007:HXDecomp, KV:2009:HCurl, AFW:2000:H(div), Toselli:2000:OSHCurl, Toselli:2004:FETIDPHCurl, Calvo:2015:OSHCurl, Calvo:2020:OSHCurl, HZ:2003:DDMaxwell, DW:2015:BDDCHCurl, Zampini:2017:BDDCHCurl, LWXZ:2007:RSCNSS}.
	
	It is well-known that the traditional smoothers for solving the scalar elliptic problems do not work well for vector field problems related to $\HCurl$ and $\HDiv$; see \cite{CGP:1993:MLMFEM}. Hence, a special smoothing technique is necessary for vector field problems. Function space splitting methods based on Helmholtz type decompositions pioneered by Hiptmair \cite{Hiptmair:1999:MGHCurl, Hiptmair:1997:MGHDiv} and Hiptmair and Xu \cite{HX:2007:HXDecomp} have been considered. In \cite{AFW:1997:H(div), AFW:2000:H(div), AFW:1998:H(DIV)}, an overlapping type domain decomposition preconditioner has been applied. We also note that the author and Brenner considered smoothers associated with nonoverlapping type domain decomposition methods for $\HDiv$ problems in \cite{BO:2018:MGHdiv, BO:2018:MGHdivNE}.  
	
	In this paper, we propose a $V$--cycle multigrid method with nonoverlapping domain decomposition smoothers and mainly consider the numerical study that is not covered by theories in \cite{Oh:2022:MGHCurl}. In \cite{Oh:2022:MGHCurl}, the author provided the convergence analysis with the assumptions, i.e., convex domain and constant material parameters in \eqref{eq:ModelProblem}. We test our method with less strict conditions, e.g., jump coefficients for $\alpha$ and $\beta$, nonconvex domain. We note that our multigrid method is an $\HCurl$ counterpart of the method in \cite{BO:2018:MGHdiv, BO:2018:MGHdivNE} and a nonoverlapping alternative of the method in \cite{AFW:2000:H(div)}, which requires less computational costs when applying the smoother.
	
	The remainder of this paper is structured as follows. In Section~\ref{sec:ND}, the standard way to discretize our model problem using the hexahedral N\'{e}d\'{e}lec element is introduced. We next present our $V$--cycle multigrid algorithm in Section~\ref{sec:MGM}. Finally, we provide numerical experiments in Section~\ref{sec:Numerics}.
	
	\section{The Discrete Problem}\label{sec:ND}
	We first consider a triangulation $\cT_h$ of $\Omega$ into hexahedral elements. The lowest order hexahedral N\'{e}d\'{e}lec element \cite{Monk:2003:BookMaxwell, Nedlec:1980:FEM} has the following form:
	\begin{equation*}
	\begin{bmatrix}
	a_1+a_2y+a_3z+a_4 yz\\
	b_1+b_2z+b_3x+b_4 zx\\
	c_1+c_2x+c_3y+c_4 xy
	\end{bmatrix}
	\end{equation*}
	on each hexahedral mesh, where the $a_i$'s, $b_i$'s and $c_i$'s are constants. We note that the twelve degrees of freedom can be completely recovered by the average tangential component on each edge of the element. Using the finite elements, we obtain the following discretized problem for (\ref{eq:ModelProblem}): Find $\bu_h \in W_h$ such that 
	\begin{equation}\label{eq:FEM}
		a(\bu_h,\bv_h)=(\bm{f}, \bv_h)\qquad\forall\,\bv_h\in W_h,
	\end{equation}
	where $W_h$ is the N\'{e}d\'{e}lec finite element space of the lowest order. 
	
	We define the operator $A : W_h \rightarrow W_h'$ in the following way:
	\begin{equation}\label{eq:ADef}
	\langle A\bv_h,\bw_h\rangle=a(\bv_h,\bw_h)\qquad\forall\,\bv_h,\bw_h\in W_h.
	\end{equation}
	Here, $\langle\cdot,\cdot,\rangle$ is the canonical bilinear form on $W_h'\times W_h$. We then have the following discrete problem:
	\begin{equation}\label{eq:discrete system}
	A\bu_h = f_h,
	\end{equation}
	where $f_h \in W_h'$ defined by 
	 \begin{equation*}
	\langle f_h,\bv_h\rangle=(\bm{f},\bv_h)\qquad\forall\,\bv_h\in W_h.
	\end{equation*}
	\section{$V$-Cycle Multigrid Method}\label{sec:MGM}
	In a multigrid setting, we construct a nested sequence of triangulations, $\cT_0, \cT_1, \cdots, $ starting with the initial triangulation $\cT_0$ consisting of few hexahedral elements. We assume that $\cT_{k}$ is obtained by uniform subdivision from $\cT_{k-1}$. We then define $W_k$ which is the lowest order N\'{e}d\'{e}lec space related to the $k^{\rm th}$ level mesh and the corresponding discrete problem: Find $\bu_k \in W_k$ such that
	\begin{equation}\label{eq:discrete system k}
	A_k \bu_k = f_k.
	\end{equation}
	Here, $f_k$ is defined by 
	\begin{equation*}
	\langle f_k,\bv_k\rangle=(\bm{f},\bv_k)\qquad\forall\,\bv_k\in W_k.
	\end{equation*}
	
	In order for solving the discrete problem (\ref{eq:discrete system k}) using multigrid methods, we need two essential ingredients, intergrid transfer operators and smoothers. We note that everything else can be constructed in a standard way.
	
	Since we are dealing with the nested finite element spaces, the natural injection can be used as the coarse-to-fine operator $I_{k-1}^k:W_{k-1}\longrightarrow W_k$. The fine-to-coarse operator $I_k^{k-1}:W_k'\longrightarrow W_{k-1}'$
	is  then defined by
	\begin{equation}\label{eq:F2C}
	\langle I_k^{k-1}r,\bv_{k-1}\rangle=\langle r,I_{k-1}^k\bv_{k-1}\rangle
	\qquad\forall\,r\in W_k',\,\bv_{k-1}\in W_{k-1}.
	\end{equation}
	
	We now focus on the missing piece, smoother. As we are interested in designing nonoverlapping type domain decomposition smoothers, we borrow the standard two level domain decomposition framework, i.e., the coarse level and the fine level that are associated with $\cT_{k-1}$ and $\cT_{k}$, respectively.
	
	Before we construct the smoothers, we set up notations for the geometric substructures. Let $\cE_{k-1}$, $\cF_{k-1}$, and $\cV_{k-1}$ be the sets of interior edges, interior faces, and interior vertices of the triangulation $\cT_{k-1}$, respectively. 
	
	We first consider the interior space. Given any coarse element $T \in \cT_{k-1}$, let us define the subspace $W_{k}^T$ by
	 \begin{equation}\label{eq:SubdomainSpace}
	W_k^T=\{\bv\in W_k:\,\bv=\bm{0}\;\text{on}\; \Omega\setminus T\}.
	\end{equation} 
	 Let $J_T$ denote the natural injection from $W_k^T$ into $W_k$. The operator $A_T:W_k^T\rightarrow \left(W_k^T\right)'$ is constructed by
	\begin{equation}\label{eq:ATDef}
		\langle A_T\bv,\bw\rangle= a(\bv,\bw)\qquad\forall\,\bv,\bw\in W_k^T.
	\end{equation}
	For a coarse edge $E \in \cE_{k-1}$, there are four coarse elements, $T_E^i, i = 1, 2, 3, 4$, in $\cT_{k-1}$ and four coarse faces, $F_E^i, i = 1, 2, 3, 4$, in $\cF_{k-1}$, that are sharing $E$. The edge space $W_{k}^E$ of $W_k$ is defined as follow:
	\begin{equation}
		\begin{aligned}
			W_k^E = & \left\{ \bv \in W_k : \bv = \bm{0} \;\text{ on }\; \Omega \setminus \left(\left(\cup_{i=1}^4 T_E^i\right) \bigcup \left(\cup_{j=1}^4 F_E^j\right) \bigcup E \right), \right.\\ 
			& \left. \hspace{10pt} \mbox{ and } a(\bv, \bw)  = 0 \quad \forall \, \bw \in 
			\left( W_k^{T_E^1} + W_k^{T_E^2} + W_k^{T_E^3} + W_k^{T_E^4} \right) \right\}.
		\end{aligned}
	\end{equation}
	 Let $J_E:W_k^E\rightarrow W_k$ be the natural injection and the operator $A_E : W_k^E \rightarrow \left(W_k^E\right)'$ be defined by
	\begin{equation}\label{eq:AEDef}
		\langle A_E\bv,\bw\rangle= a(\bv,\bw)\qquad\forall\,\bv,\bw\in W_k^E.
	\end{equation}
	Finally, we define the vertex space $W_{k}^P$ of $W_k$. For each coarse vertex $P \in \cV_{k-1}$, there are eight elements, $T_P^i, i = 1, \cdots, 8$, in $\cT_{k-1}$, twelve faces, $F_P^j, j = 1, \cdots, 12$, in $\cF_{k-1}$, and six edges, $E_P^l, l = 1, \cdots, 6$, in $\cE_{k-1}$, that have the point $P$ in common. We define the vertex space $W_{k}^P$ by
	\begin{equation}
		\begin{aligned}
			W_k^P = & \left\{ \bv \in W_k : \bv = \bm{0} \;\text{ on }\; \Omega \setminus \left(\left(\cup_{i=1}^8 T_P^i\right) \bigcup \left(\cup_{j=1}^{12} F_P^j\right) \bigcup \left(\cup_{l=1}^6 E_P^l\right) \right), \right.\\ 
			& \left. \hspace{10pt} \mbox{ and } a(\bv, \bw)  = 0 \quad \forall \, \bw \in 
			\left( \sum_{i=1}^8 W_k^{T_P^i}\right) \right\}.
		\end{aligned}
	\end{equation}
	The natural injection $J_P:W_k^P\rightarrow W_k$ and the operator $A_P$ are obtained by a similar way to $J_E$ and $A_E$, respectively.
	
	We now define two smoothers, the edge-based and the vertex-based preconditioners. The edge-based smoother $M_{E, k}^{-1}$ is constructed as follow:
	\begin{equation}\label{eq:Preconditioner E}
		M_{E,k}^{-1}=\eta_E\left(\sum_{T\in\cT_{k-1}}J_T A_T^{-1} J_T^t+\sum_{E\in\cE_{k-1}}J_E A_E^{-1} J_E^t\right).
	\end{equation}
	Similarly, the vertex-based smoother $M_{P, k}^{-1}$ is obtained by
	\begin{equation}\label{eq:Preconditioner P}
		M_{P,k}^{-1}=\eta_P\left(\sum_{T\in\cT_{k-1}}J_T A_T^{-1} J_T^t+\sum_{P\in\cV_{k-1}}J_P A_P^{-1} J_P^t\right).
	\end{equation}
	Here, $\eta_E$ and $\eta_P$ are damping factors and $J_T^t:W_k'\rightarrow \left(W_k^T\right)'$, $J_E^t:W_k'\rightarrow \left(W_k^E\right)'$, and $J_P^t:W_k'\rightarrow \left(W_k^P\right)'$ are the transposes of $J_T$, $J_E$, and $J_P$, respectively. We can decide the damping factors such that
	\begin{equation}\label{eq:spectral condition}
	\rho\left(M_{E,k}^{-1}A_k\right) \le 1 \mbox{ and } \rho\left(M_{P,k}^{-1} A_k\right) \le 1,
	\end{equation}
	where $\rho\left(M_{E,k}^{-1}A_k\right)$ and $\rho\left(M_{P,k}^{-1}A_k\right)$ are the spectral radii of $M_{E,k}^{-1}A_k$ and $M_{P,k}^{-1}A_k$, respectively. We note that the conditions in \eqref{eq:spectral condition} are satisfied if $\eta_E \le 1/12$ and $\eta_P \le 1/8$.
	
	Putting all together, we can completely determine the multigrid framework.
	The output $MGV\left(k, g, \bz_0, m\right)$ of the $k^{\rm th}$ level (symmetric) multigrid $V$-cycle algorithm
	for $A_k \bz = g$
	, with initial guess $\bz_0\in W_k$ and $m$ smoothing steps,
	is defined as follows:
	\par
	For $k=0$, the result is obtained by using a direct solver:
	\begin{equation*}
	MGV\left(0, g, \bz_0, m\right) = A_0^{-1}g.\nonumber
	\end{equation*}
	\par\noindent
	For $k\ge1$, we set
	\begin{align*}
	\bz_l &= \bz_{l-1} + M_k^{-1} \left( g - A_k \bz_{l-1} \right) \qquad \text{for}\; l = 1, \cdots, m,\nonumber\\
	\overline{g} &= I_{k}^{k-1} \left(g - A_k \bz_{m} \right),\nonumber\\
	\bz_{m+1} &= \bz_{m} + I_{k-1}^k MGV\left(k-1, \overline{g}, 0, m\right), \nonumber\\
	\bz_l &= \bz_{l-1} + M_k^{-1} \left( g - A_k \bz_{l-1} \right)\qquad\text{for}\; l = m + 2, \cdots, 2m + 1.\nonumber
	\end{align*}
	The output of $MGV\left(k, g, \bz_0, m\right)$ is $\bz_{2m + 1}$.
	The smoother $M_k$ is either $M_{E, k}$ or $M_{P, k}$ defined earlier in the current section.
	
	In order to check the efficiency of the algorithm, we consider the following error propagation operator:
	\begin{equation}\label{eq:error propagation operator}
		\mathbb{E}_k^m ( \bz - \bz_0) = \bz - MGV(k, g, \bz_0, m).
	\end{equation}
	The operator $\mathbb{E}_k^m$ is affiliated with the error after one $k^{\rm th}$ multigrid sweep with $m$ smoothing steps. For more detail, see \cite[Chapter 6]{Brenner:2008:book}.
	
	\section{Numerical Results}\label{sec:Numerics}
	We note that a part of the finite element discretizations has been implemented with the MFEM library \cite{mfem, mfem-web}. The codes used for the experiments are available at the repository \url{https://github.com/duksoon-open/MG_ND}.
	
	\subsection{Jump Coefficients}\label{subsec:jump coeff}
	\begin{figure}[!ht]
		\centering
		\includegraphics[scale=0.10]{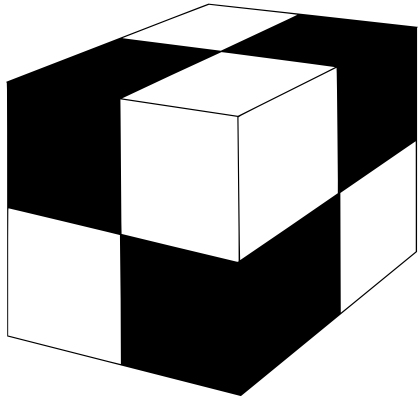}
		\caption{Checkerboard distribution of the coefficients}
		\label{fig:checkerboard}
	\end{figure}
	
	The first experiment is for the cube $\Omega = (-1, 1)^3$. The domain $\Omega$ is decomposed into eight identical cubical subdomains and set the subregions as the initial triangulation $\cT_0$. We assume that the coefficients $\alpha$ and $\beta$ are constants in each subdomain in a checkerboard pattern; see Fig.~\ref{fig:checkerboard}. We apply the multigrid algorithms, edge-based method and vertex-based method, introduced in Section~\ref{sec:MGM}. We report the contraction numbers by calculating the largest eigenvalue of the operators $\mathbb{E}_k^m$ in \eqref{eq:error propagation operator}, for $k=1,\cdots,4$ and for $m=1,\cdots,5$. The results are presented in Table~\ref{table:Jump Coefficient Edge} and Table~\ref{table:Jump Coefficient Vertex}. We see that the multigrid methods provide contraction and are robust to the jump between the interface of the initial mesh.
	
	\begin{table}[!ht]
		\small
		\caption{Contraction numbers of the $V$-cycle edge-based methods. $\alpha_b$ and $\beta_b$ for the black subregions and $\alpha_w$ and $\beta_w$ for the white subregions as indicated in a checkerboard pattern as in Fig.~\ref{fig:checkerboard}}
		\centering
		\begin{tabular}{c|ccccc}
			& $m=1$  & $m=2$  & $m=3$  & $m=4$  & $m=5$   \\ \hline
			& \multicolumn{5}{ c }{$\alpha_b = 0.01, \beta_b = 1, \alpha_w = 1, \beta_w = 1$}  \\ \hline
			$k=1$ &0.905&0.827&0.762&0.709&0.663\\ \hline
			$k=2$ &0.940&0.908&0.872&0.841&0.811\\ \hline
			$k=3$ &0.967&0.952&0.935&0.917&0.902\\ \hline
			$k=4$ &0.981&0.970&0.960&0.955&0.942\\ \hline \hline
			& \multicolumn{5}{ c }{$\alpha_b = 0.1, \beta_b = 1, \alpha_w = 1, \beta_w = 1$}  \\ \hline
			$k=1$ &0.905&0.827&0.763&0.710&0.666\\ \hline
			$k=2$ &0.941&0.910&0.875&0.844&0.807\\ \hline
			$k=3$ &0.967&0.954&0.937&0.92&0.905\\ \hline
			$k=4$ &0.980&0.971&0.961&0.956&0.945\\ \hline \hline
			& \multicolumn{5}{ c }{$\alpha_b = 1, \beta_b = 1, \alpha_w = 1, \beta_w = 1$}  \\ \hline
			$k=1$ &0.907&0.831&0.769&0.719&0.677\\ \hline
			$k=2$ &0.944&0.917&0.885&0.858&0.830\\ \hline
			$k=3$ &0.970&0.959&0.944&0.930&0.917\\ \hline
			$k=4$ &0.981&0.972&0.965&0.963&0.956\\ \hline \hline
			& \multicolumn{5}{ c }{$\alpha_b = 10, \beta_b = 1, \alpha_w = 1, \beta_w = 1$}  \\ \hline
			$k=1$ &0.909&0.836&0.777&0.729&0.690\\ \hline
			$k=2$ &0.948&0.924&0.896&0.872&0.853\\ \hline
			$k=3$ &0.972&0.965&0.952&0.941&0.931\\ \hline
			$k=4$ &0.982&0.974&0.971&0.969&0.966\\ \hline \hline
			& \multicolumn{5}{ c }{$\alpha_b = 100, \beta_b = 1, \alpha_w = 1, \beta_w = 1$}  \\ \hline
			$k=1$ &0.910&0.837&0.778&0.731&0.693\\ \hline
			$k=2$ &0.949&0.926&0.898&0.875&0.857\\ \hline
			$k=3$ &0.973&0.966&0.954&0.943&0.934\\ \hline
			$k=4$ &0.982&0.975&0.972&0.970&0.968\\ \hline 
		\end{tabular}
		\label{table:Jump Coefficient Edge}
	\end{table}
	
	\begin{table}[!ht]
		\small
		\caption{Contraction numbers of the $V$-cycle vertex-based methods. $\alpha_b$ and $\beta_b$ for the black subregions and $\alpha_w$ and $\beta_w$ for the white subregions as indicated in a checkerboard pattern as in Fig.~\ref{fig:checkerboard}}
		\centering
		\begin{tabular}{c|ccccc}
			& $m=1$  & $m=2$  & $m=3$  & $m=4$  & $m=5$   \\ \hline
			& \multicolumn{5}{ c }{$\alpha_b = 0.01, \beta_b = 1, \alpha_w = 1, \beta_w = 1$}  \\ \hline
			$k=1$ &0.790&0.624&0.493&0.390&0.308\\ \hline
			$k=2$ &0.792&0.627&0.495&0.393&0.312\\ \hline
			$k=3$ &0.791&0.625&0.494&0.391&0.310\\ \hline
			$k=4$ &0.791&0.626&0.495&0.392&0.317\\ \hline \hline
			& \multicolumn{5}{ c }{$\alpha_b = 0.1, \beta_b = 1, \alpha_w = 1, \beta_w = 1$}  \\ \hline
			$k=1$ &0.790&0.624&0.493&0.390&0.308\\ \hline
			$k=2$ &0.791&0.626&0.494&0.392&0.310\\ \hline
			$k=3$ &0.791&0.626&0.495&0.392&0.310\\ \hline
			$k=4$ &0.791&0.626&0.495&0.392&0.311\\ \hline \hline
			& \multicolumn{5}{ c }{$\alpha_b = 1, \beta_b = 1, \alpha_w = 1, \beta_w = 1$}  \\ \hline
			$k=1$ &0.790&0.624&0.493&0.390&0.308\\ \hline
			$k=2$ &0.791&0.626&0.495&0.392&0.310\\ \hline
			$k=3$ &0.791&0.626&0.495&0.392&0.311\\ \hline
			$k=4$ &0.791&0.626&0.495&0.392&0.311\\ \hline \hline
			& \multicolumn{5}{ c }{$\alpha_b = 10, \beta_b = 1, \alpha_w = 1, \beta_w = 1$}  \\ \hline
			$k=1$ &0.790&0.624&0.493&0.390&0.308\\ \hline
			$k=2$ &0.791&0.626&0.495&0.392&0.311\\ \hline
			$k=3$ &0.791&0.626&0.495&0.392&0.311\\ \hline
			$k=4$ &0.791&0.626&0.495&0.392&0.311\\ \hline \hline
			& \multicolumn{5}{ c }{$\alpha_b = 100, \beta_b = 1, \alpha_w = 1, \beta_w = 1$}  \\ \hline
			$k=1$ &0.790&0.624&0.493&0.390&0.308\\ \hline
			$k=2$ &0.791&0.626&0.495&0.392&0.311\\ \hline
			$k=3$ &0.791&0.626&0.495&0.393&0.318\\ \hline
			$k=4$ &0.791&0.632&0.524&0.483&0.429\\ \hline
		\end{tabular}
		\label{table:Jump Coefficient Vertex}
	\end{table}
	
	\subsection{Nonconvex Domain}
\begin{figure}[!ht]
		\centering
		\includegraphics[scale=0.10]{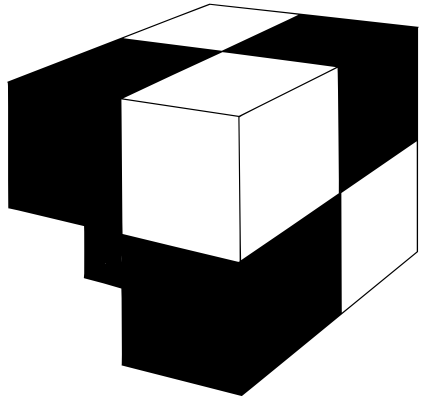}
		\caption{Checkerboard distribution of the coefficients for the nonconvex domain}
		\label{fig:checkerboard_nc}
	\end{figure}

	In the second set of numerical tests, we consider the domain $\Omega = (-1, 1)^3 \backslash (-1, 0)^3$, which is nonconvex. We begin with the initial mesh $\cT_{0}$ which consists of seven identical cubes as in Figure~\ref{fig:checkerboard_nc} and inductively define the $k^{\rm th}$ level mesh $\cT_k$ by a uniform subdivision.
	Other general settings are similar to those of Section~\ref{subsec:jump coeff}. 
	We apply our multigrid algorithm for the problem (\ref{eq:ModelProblem}) on the domain $\Omega$ and report the contraction numbers computed in the same manner with the first experiment. As we see the results in the Table~\ref{table:Nonconvex Jump Coefficient Edge} and Table~\ref{table:Nonconvex Jump Coefficient Vertex}, the uniform convergences and robustness are observed except for the problem with $k=1$ using the vertex-based smoother.
	
	\begin{table}[!ht]
		\small
		\caption{Contraction numbers of the $V$-cycle edge-based methods for the nonconvex domain. $\alpha_b$ and $\beta_b$ for the black subregions and $\alpha_w$ and $\beta_w$ for the white subregions as indicated in a checkerboard pattern as in Fig.~\ref{fig:checkerboard_nc}}
		\centering
		\begin{tabular}{c|ccccc}
			& $m=1$  & $m=2$  & $m=3$  & $m=4$  & $m=5$   \\ \hline
			& \multicolumn{5}{ c }{$\alpha_b = 0.01, \beta_b = 1, \alpha_w = 1, \beta_w = 1$}  \\ \hline
			$k=1$ &0.835&0.702&0.576&0.471&0.420\\ \hline
			$k=2$ &0.940&0.881&0.828&0.785&0.749\\ \hline
			$k=3$ &0.967&0.940&0.918&0.892&0.869\\ \hline
			$k=4$ &0.982&0.969&0.957&0.948&0.934\\ \hline \hline
			& \multicolumn{5}{ c }{$\alpha_b = 0.1, \beta_b = 1, \alpha_w = 1, \beta_w = 1$}  \\ \hline
			$k=1$ &0.823&0.682&0.542&0.436&0.390\\ \hline
			$k=2$ &0.943&0.887&0.832&0.796&0.757\\ \hline
			$k=3$ &0.966&0.941&0.918&0.892&0.870\\ \hline
			$k=4$ &0.983&0.970&0.958&0.948&0.936\\ \hline \hline
			& \multicolumn{5}{ c }{$\alpha_b = 1, \beta_b = 1, \alpha_w = 1, \beta_w = 1$}  \\ \hline
			$k=1$ &0.799&0.627&0.511&0.419&0.338\\ \hline
			$k=2$ &0.945&0.894&0.844&0.802&0.768\\ \hline
			$k=3$ &0.970&0.942&0.920&0.896&0.876\\ \hline
			$k=4$ &0.984&0.974&0.961&0.950&0.941\\ \hline \hline
			& \multicolumn{5}{ c }{$\alpha_b = 10, \beta_b = 1, \alpha_w = 1, \beta_w = 1$}  \\ \hline
			$k=1$ &0.802&0.646&0.523&0.420&0.348\\ \hline
			$k=2$ &0.946&0.896&0.851&0.804&0.773\\ \hline
			$k=3$ &0.972&0.943&0.924&0.902&0.881\\ \hline
			$k=4$ &0.986&0.976&0.964&0.954&0.945\\ \hline \hline
			& \multicolumn{5}{ c }{$\alpha_b = 100, \beta_b = 1, \alpha_w = 1, \beta_w = 1$}  \\ \hline
			$k=1$ &0.804&0.651&0.529&0.424&0.354\\ \hline
			$k=2$ &0.947&0.896&0.851&0.805&0.774\\ \hline
			$k=3$ &0.972&0.944&0.924&0.903&0.882\\ \hline
			$k=4$ &0.986&0.976&0.965&0.955&0.946\\ \hline
		\end{tabular}
		\label{table:Nonconvex Jump Coefficient Edge}
	\end{table}
	
	\begin{table}[!ht]
		\small
		\caption{Contraction numbers of the $V$-cycle vertex-based methods for the nonconvex domain. $\alpha_b$ and $\beta_b$ for the black subregions and $\alpha_w$ and $\beta_w$ for the white subregions as indicated in a checkerboard pattern as in Fig.~\ref{fig:checkerboard_nc}}
		\centering
		\begin{tabular}{c|ccccc}
			& $m=1$  & $m=2$  & $m=3$  & $m=4$  & $m=5$ \\ \hline
			& \multicolumn{5}{ c }{$\alpha_b = 0.01, \beta_b = 1, \alpha_w = 1, \beta_w = 1$}  \\ \hline
			$k=1$ &$>1$&$>1$&$>1$&$>1$&$>1$\\ \hline
			$k=2$ &0.912&0.835&0.765&0.683&0.648\\ \hline
			$k=3$ &0.905&0.825&0.758&0.689&0.622\\ \hline
			$k=4$ &0.912&0.834&0.765&0.697&0.626\\ \hline \hline
			& \multicolumn{5}{ c }{$\alpha_b = 0.1, \beta_b = 1, \alpha_w = 1, \beta_w = 1$}  \\ \hline
			$k=1$ &$>1$&$>1$&$>1$&$>1$&$>1$\\ \hline
			$k=2$ &0.896&0.808&0.726&0.622&0.598\\ \hline
			$k=3$ &0.878&0.789&0.710&0.629&0.535\\ \hline
			$k=4$ &0.890&0.792&0.703&0.577&0.518\\ \hline \hline
			& \multicolumn{5}{ c }{$\alpha_b = 1, \beta_b = 1, \alpha_w = 1, \beta_w = 1$}  \\ \hline
			$k=1$ &$>1$&$>1$&$>1$&$>1$&$>1$\\ \hline
			$k=2$ &0.869&0.758&0.660&0.576&0.504\\ \hline
			$k=3$ &0.837&0.705&0.593&0.497&0.425\\ \hline
			$k=4$ &0.827&0.689&0.570&0.479&0.407\\ \hline \hline
			& \multicolumn{5}{ c }{$\alpha_b = 10, \beta_b = 1, \alpha_w = 1, \beta_w = 1$}  \\ \hline
			$k=1$ &$>1$&$>1$&$>1$&$>1$&$>1$\\ \hline
			$k=2$ &0.863&0.766&0.666&0.591&0.503\\ \hline
			$k=3$ &0.844&0.725&0.613&0.495&0.453\\ \hline
			$k=4$ &0.829&0.706&0.574&0.498&0.420\\ \hline \hline
			& \multicolumn{5}{ c }{$\alpha_b = 100, \beta_b = 1, \alpha_w = 1, \beta_w = 1$}  \\ \hline
			$k=1$ &$>1$&$>1$&$>1$&$>1$&$>1$\\ \hline
			$k=2$ &0.865&0.772&0.675&0.602&0.515\\ \hline
			$k=3$ &0.849&0.735&0.626&0.509&0.471\\ \hline
			$k=4$ &0.833&0.717&0.590&0.519&0.447\\ \hline
		\end{tabular}
		\label{table:Nonconvex Jump Coefficient Vertex}
	\end{table}
	
	\section*{Acknowledgement}
	This work was supported by the National Research Foundation of Korea (NRF) grant funded by the Korea government (MSIT) (No. 2020R1F1A1A01072168).


\begin{thebibliography}{10}

		\bibitem{mfem}
		R.~Anderson, J.~Andrej, A.~Barker, J.~Bramwell, J.-S. Camier, J.~C.~V. Dobrev,
		  Y.~Dudouit, A.~Fisher, T.~Kolev, W.~Pazner, M.~Stowell, V.~Tomov,
		  I.~Akkerman, J.~Dahm, D.~Medina, and S.~Zampini.
		\newblock {MFEM}: A modular finite element library.
		\newblock {\em Computers \& Mathematics with Applications}, 81:42--74, 2021.
		
		\bibitem{AFW:1997:H(div)}
		D.~N. Arnold, R.~S. Falk, and R.~Winther.
		\newblock Preconditioning in {$H({\rm div})$} and applications.
		\newblock {\em Math. Comp.}, 66(219):957--984, 1997.
		
		\bibitem{AFW:1998:H(DIV)}
		D.~N. Arnold, R.~S. Falk, and R.~Winther.
		\newblock Multigrid preconditioning in {$H({\rm div})$} on non-convex polygons.
		\newblock {\em Comput. Appl. Math.}, 17(3):303--315, 1998.
		
		\bibitem{AFW:2000:H(div)}
		D.~N. Arnold, R.~S. Falk, and R.~Winther.
		\newblock Multigrid in {$H({\rm div})$} and {$H({\rm curl})$}.
		\newblock {\em Numer. Math.}, 85(2):197--217, 2000.
		
		\bibitem{BO:2018:MGHdiv}
		S.~C. Brenner and D.-S. Oh.
		\newblock Multigrid methods for {$H({\rm div})$} in three dimensions with
		  nonoverlapping domain decomposition smoothers.
		\newblock {\em Numer. Linear Algebra Appl.}, 25(5):e2191, 14, 2018.
		
		\bibitem{BO:2018:MGHdivNE}
		S.~C. Brenner and D.-S. Oh.
		\newblock A smoother based on nonoverlapping domain decomposition methods for
		  {$H({\rm div})$} problems: a numerical study.
		\newblock In {\em Domain decomposition methods in science and engineering
		  {XXIV}}, volume 125 of {\em Lect. Notes Comput. Sci. Eng.}, pages 523--531.
		  Springer, Cham, 2018.
		
		\bibitem{Brenner:2008:book}
		S.~C. Brenner and L.~R. Scott.
		\newblock {\em The mathematical theory of finite element methods}, volume~15 of
		  {\em Texts in Applied Mathematics}.
		\newblock Springer, New York, third edition, 2008.
		
		\bibitem{CGP:1993:MLMFEM}
		Z.~Q. Cai, C.~I. Goldstein, and J.~E. Pasciak.
		\newblock Multilevel iteration for mixed finite element systems with penalty.
		\newblock {\em SIAM J. Sci. Comput.}, 14(5):1072--1088, 1993.
		
		\bibitem{Calvo:2015:OSHCurl}
		J.~G. Calvo.
		\newblock A two-level overlapping {S}chwarz method for {$H(\rm curl)$} in two
		  dimensions with irregular subdomains.
		\newblock {\em Electron. Trans. Numer. Anal.}, 44:497--521, 2015.
		
		\bibitem{Calvo:2020:OSHCurl}
		J.~G. Calvo.
		\newblock A new coarse space for overlapping {S}chwarz algorithms for {$H(\rm
		  curl)$} problems in three dimensions with irregular subdomains.
		\newblock {\em Numer. Algorithms}, 83(3):885--899, 2020.
		
		\bibitem{DW:2015:BDDCHCurl}
		C.~R. Dohrmann and O.~B. Widlund.
		\newblock A {BDDC} algorithm with deluxe scaling for three-dimensional {$H({\bf
		  curl})$} problems.
		\newblock {\em Comm. Pure Appl. Math.}, 69(4):745--770, 2016.
		
		\bibitem{GR:1986:NS}
		V.~Girault and P.-A. Raviart.
		\newblock {\em Finite element methods for {N}avier-{S}tokes equations},
		  volume~5 of {\em Springer Series in Computational Mathematics}.
		\newblock Springer-Verlag, Berlin, 1986.
		\newblock Theory and algorithms.
		
		\bibitem{Hiptmair:1997:MGHDiv}
		R.~Hiptmair.
		\newblock Multigrid method for {$H({\rm div})$} in three dimensions.
		\newblock {\em Electron. Trans. Numer. Anal.}, 6(Dec.):133--152, 1997.
		\newblock Special issue on multilevel methods (Copper Mountain, CO, 1997).
		
		\bibitem{Hiptmair:1999:MGHCurl}
		R.~Hiptmair.
		\newblock Multigrid method for {M}axwell's equations.
		\newblock {\em SIAM J. Numer. Anal.}, 36(1):204--225, 1999.
		
		\bibitem{HX:2007:HXDecomp}
		R.~Hiptmair and J.~Xu.
		\newblock Nodal auxiliary space preconditioning in {${\bf H}({\bf curl})$} and
		  {${\bf H}({\rm div})$} spaces.
		\newblock {\em SIAM J. Numer. Anal.}, 45(6):2483--2509, 2007.
		
		\bibitem{HZ:2003:DDMaxwell}
		Q.~Hu and J.~Zou.
		\newblock A nonoverlapping domain decomposition method for {M}axwell's
		  equations in three dimensions.
		\newblock {\em SIAM J. Numer. Anal.}, 41(5):1682--1708, 2003.
		
		\bibitem{KV:2009:HCurl}
		T.~V. Kolev and P.~S. Vassilevski.
		\newblock Parallel auxiliary space {AMG} for {$H({\rm curl})$} problems.
		\newblock {\em J. Comput. Math.}, 27(5):604--623, 2009.
		
		\bibitem{LWXZ:2007:RSCNSS}
		Y.-J. Lee, J.~Wu, J.~Xu, and L.~Zikatanov.
		\newblock Robust subspace correction methods for nearly singular systems.
		\newblock {\em Math. Models Methods Appl. Sci.}, 17(11):1937--1963, 2007.
		
		\bibitem{mfem-web}
		{MFEM}: Modular finite element methods {[Software]}.
		\newblock \url{mfem.org}.
		
		\bibitem{Monk:2003:BookMaxwell}
		P.~Monk.
		\newblock {\em Finite element methods for {M}axwell's equations}.
		\newblock Numerical Mathematics and Scientific Computation. Oxford University
		  Press, New York, 2003.
		
		\bibitem{Monk:1991:Maxwell}
		P.~B. Monk.
		\newblock A mixed method for approximating {M}axwell's equations.
		\newblock {\em SIAM J. Numer. Anal.}, 28(6):1610--1634, 1991.
		
		\bibitem{Nedlec:1980:FEM}
		J.-C. N\'{e}d\'{e}lec.
		\newblock Mixed finite elements in {${\bf R}^{3}$}.
		\newblock {\em Numer. Math.}, 35(3):315--341, 1980.
		
		\bibitem{Oh:2022:MGHCurl}
		D.-S. Oh.
		\newblock Multigrid methods for 3{$D$} ${H}(\mathbf{curl})$ problems with
		  nonoverlapping domain decomposition smoothers, 2022.
		\newblock submitted, arXiv:2205.05840.
		
		\bibitem{Toselli:2000:OSHCurl}
		A.~Toselli.
		\newblock Overlapping {S}chwarz methods for {M}axwell's equations in three
		  dimensions.
		\newblock {\em Numer. Math.}, 86(4):733--752, 2000.
		
		\bibitem{Toselli:2004:FETIDPHCurl}
		A.~Toselli.
		\newblock Domain decomposition methods of dual-primal {FETI} type for edge
		  element approximations in three dimensions.
		\newblock {\em C. R. Math. Acad. Sci. Paris}, 339(9):673--678, 2004.
		
		\bibitem{Zampini:2017:BDDCHCurl}
		S.~Zampini.
		\newblock Adaptive {BDDC} deluxe methods for {$\rm H(curl)$}.
		\newblock In {\em Domain decomposition methods in science and engineering
		  {XXIII}}, volume 116 of {\em Lect. Notes Comput. Sci. Eng.}, pages 285--292.
		  Springer, Cham, 2017.
		
		\end{thebibliography}
\end{document}